# PREHOMOGENEOUS VECTOR SPACES AND ERGODIC THEORY III

Akihiko Yukie[1]

Oklahoma State University

**Contents**

**Introduction**

This is part three in this series of papers. Throughout this paper, $k$ is a field of characteristic zero. In parts one and two [25], [23], we consider cases (2), (5), (6), (7) in the Sato–Kimura classification [19]. In this part, we consider the prehomogeneous vector space $G = \text{GL}(5) \times \text{GL}(3)$, $V = \wedge^2 k^5 \otimes k^3$. This is rather an isolated case and that's why we consider this single case separately here. However, as far as the relation between the theory of prehomogeneous vector spaces and ergodic theory is concerned, this may be the most interesting case besides the situation of the Oppenheim conjecture $G = \text{GL}(n)$, $V = \text{Sym}^2(k^n)^*$, because it produces a family of (irrational) cubic forms in five variables whose values at integer points are dense in $\mathbb{R}$. This may be the family of cubic forms in the lowest number of variables we can achieve by the theory of prehomogeneous vector spaces. We consider most of the remaining irreducible split cases in part four.

Oppenheim conjectured in [12] that if $Q(x)$ is a real non-degenerate indefinite quadratic form in $n \geq 5$ variables such that the ratio of at least one pair of coefficients is irrational, for any $\epsilon > 0$, there exists $x \in \mathbb{Z}^n$ such that $0 < |Q(x)| < \epsilon$. Due to the result of Lewis [9], this is equivalent to say the set $\{Q(x) \mid x \in \mathbb{Z}^n\}$ is dense in $\mathbb{R}$. There were many partial results including the one by Davenport with the collaboration with others ([4], [5], [6], [2], [18]) for $n \geq 21$. It was proved in the optimum form by Margulis (see [10]) for $n \geq 3$ using ergodic theory.

We posed the question of generalizing the Oppenheim conjecture from the viewpoint of prehomogeneous vector spaces in [25]. For more detailed comments, the

[1]Partially supported by NSF grant DMS-9401391

Typeset by $\mathcal{AMS}$-TEX

reader should see the introduction of [25]. Here, we briefly state what we are going to prove.

Let $H_1 = \mathrm{SL}(5), H_2 = \mathrm{SL}(3), H = H_1 \times H_2$. It is known that $(G, V)$ is a prehomogeneous vector space (see [22], [26], [25] for the definition of prehomogeneous vector spaces). A non-constant polynomial $\Delta(x)$ on $V$ is called a relative invariant polynomial if there exists a character $\chi$ such that $\Delta(gx) = \chi(g)\Delta(x)$. Such $\Delta(x)$ exists for our case and is essentially unique. So we define $V^{\mathrm{ss}} = \{x \in V \mid \Delta(x) \neq 0\}$. For $x \in V^{\mathrm{ss}}_{\mathbb{R}}$, let $H^0_{x\mathbb{R}+}$ be the connected component of 1 in classical topology of the stabilizer $H_{x\mathbb{R}}$. We will prove that if $x \in V^{\mathrm{ss}}_{\mathbb{R}}$ is "sufficiently irrational", $H^0_{x\mathbb{R}+} H_{\mathbb{Z}}$ is dense in $H_{\mathbb{R}}$.

What Margulis did was to prove the above statement for the case $H = \mathrm{SL}(3), V = \mathrm{Sym}^2(\mathbb{R}^3)^*$. Our method is based on the following theorem due to Ratner.

**Theorem (0.1) (Ratner)** *Let $G$ be a connected Lie group and $U$ a connected subgroup of $G$ generated by unipotent elements of $G$. Then given any lattice $\Gamma \subset G$ and $x \in G/\Gamma$, there exists a connected closed subgroup $U \subset F \subset G$ such that $\overline{Ux\Gamma} = Fx\Gamma$. Moreover, $F/F \cap \Gamma$ has a finite invariant measure.*

Note that in the above theorem, the definition of a lattice contains the condition that $G/\Gamma$ has a finite volume. The first statement was called Raghunathan's topological conjecture, and the second statement was proved by Ratner in conjunction with Raghunathan's topological conjecture. Raghunathan's topological conjecture was published by Dani [3] for one dimensional unipotent groups and was generalized to groups generated by unipotent elements by Margulis [10] (before Ratner's proof). The proof for the general case was given by Ratner in a series of papers [13], [14], [15], [16]. For these, there is an excellent survey article by Ratner [17].

Note that in the above theorem, if $G$ is an algebraic group over $\mathbb{Q}$ and $\Gamma$ is an arithmetic lattice, the group $F$ becomes an algebraic group defined over $\mathbb{Q}$. For this, the reader should see Proposition (3.2) [20, pp. 321–322]. It is also proved in Proposition (3.2) [20, pp. 321–322] that the radical of $F$ is a unipotent subgroup. In [20], only one lattice is considered, but one can deduce the above statement for any lattice commensurable with the lattice in [20] by a simple argument using Ratner's theorem.

We describe our result explicitly here. For any non-zero point $x$ in a vector space, we denote the point in the corresponding projective space determined by $x$ by $[x]$. Also if $S$ is a subspace of a vector space, we denote the point in the Grassmann variety determined by $S$ by $[S]$ also. Let $V_1$ be a five dimensional vector space defined over $\mathbb{Q}$. We fix a rational basis $\{m_0, \cdots, m_4\}$ for $V_1$. Let $S$ be the subspace of $\wedge^2 V_1$ spanned by the following three elements

$$m_1 \wedge m_4 - 3m_2 \wedge m_3,$$
$$m_0 \wedge m_4 - 2m_1 \wedge m_3,$$
$$m_0 \wedge m_3 - 3m_1 \wedge m_2.$$

This space is of dimension three and therefore, determines an element $[S]$ in the Grassmann $\mathrm{Gr}(3, V_1) \cong \mathrm{Gr}(3, 10)$ of three dimensional subspaces in $V_1$. Because of the rational structure of $V_1$, $\mathrm{Gr}(3, V_1)$ has a rational structure also.

Let
$$Q(a) = a_0 a_4 - \frac{1}{4} a_1 a_3 + \frac{1}{12} a_2^2,$$
$$F(a) = 72 a_0 a_2 a_4 + 9 a_1 a_2 a_3 - 2 a_2^3 - 27 a_0 a_3^2 - 27 a_1^2 a_4$$

for $a = \sum_{i=0}^{4} a_i m_i$. If we identify $V_1$ with the space of binary quartic forms by $a \to m_0 v_1^4 + \cdots + m_4 v_2^4$, $Q, F$ correspond to quadratic and cubic SL(2)–invariant polynomials. Also $[S] \in \mathrm{Gr}(3, 10)$ is the unique fixed point of SL(2).

If $g \in \mathrm{GL}(V_1)_\mathbb{R} \cong \mathrm{GL}(5)_\mathbb{R}$, it naturally acts on $\mathrm{Gr}(3, V_1)_\mathbb{R}$, $\mathbb{P}(\mathrm{Sym}^2 V_1^*)_\mathbb{R}$, and $\mathbb{P}(\mathrm{Sym}^3 V_1^*)_\mathbb{R}$. Note that $(gQ)(a) = Q(g^{-1} a)$, $(gF)(a) = F(g^{-1} a)$. Then the following theorem is the main result of this paper.

**Theorem (0.2)** *Suppose $g[S] \notin \mathrm{Gr}(3, V_1)_\mathbb{Q}$ and $g[Q] \notin \mathbb{P}(\mathrm{Sym}^2 V_1^*)_\mathbb{Q}$. Then the set of values of the cubic polynomial $F(g^{-1} a)$ at primitive integer points in $\mathbb{Z}^5$ is dense in $\mathbb{R}$.*

In §1, we consider various identifications concerning tensor products of vector spaces. If $x \in V_\mathbb{R}^{ss}$, by Ratner's theorem, there exists a closed connected subgroup $H_{x\mathbb{R}+}^0 \subset F \subset H_\mathbb{R}$ such that $\overline{H_{x\mathbb{R}+}^0 H_\mathbb{Z}} = F H_\mathbb{Z}$. In §2, we construct equivariant maps from $V^{ss}$ to various $H$–varieties. In §3, we prove that these equivariant maps are well defined and are non-trivial. These equivariant maps correspond to families of such $F$'s with the property that if $X_F$ is the corresponding $H$–variety, $F$ has a unique fixed point in $X_F$. Part of our consideration resembles the argument in [21]. In §4, we describe the orbit space to determine when $H_{x\mathbb{R}+}^0$ is generated by unipotent elements. In §5, we classify all $F$'s as above. In §6, we prove Theorem (0.2).

## §1 Preliminaries

We are going to do a lot of computations in §3 regarding symmetric tensor products of vector spaces. We fix various normalizations for that purpose in this section.

Let $W$ be a vector space over $k$ with a basis $\{e_1, \cdots, e_n\}$. Let $W^*$ be the dual space with the dual basis $\{f_1, \cdots, f_n\}$. For $a_1, \cdots, a_d \in W$, we define

$$[a_1, \cdots, a_d]_d = \frac{1}{d!} \sum_{\sigma \in \mathfrak{S}_d} a_{\sigma(1)} \otimes \cdots \otimes a_{\sigma(d)},$$

where $\mathfrak{S}_d$ is the group of permutations of $\{1, \cdots, d\}$. We identify $\mathrm{Sym}^d W$ with the subspace of $W^{\otimes d}$ spanned by elements of the form $[a_1, \cdots, a_d]_d$. Similarly, we identify $\mathrm{Sym}^d W^*$ with a subspace of $(W^*)^{\otimes d}$. For $a_1, \cdots, a_{d_1}, a_{d_1+1}, \cdots, a_{d_1+d_2} \in W$, we define

$$[a_1, \cdots, a_{d_1}]_{d_1} [a_{d_1+1}, \cdots, a_{d_1+d_2}]_{d_2} = [a_1, \cdots, a_{d_1+d_2}]_{d_1+d_2}.$$

By this product, $\oplus \mathrm{Sym}^* W$ becomes an associative algebra. Since $[a_1, \cdots, a_d]_d = a_1 \cdots a_d$, we use this usual notation of product from now on.

Since $(W^*)^{\otimes d}$ and $W^{\otimes d}$ are dual spaces of each other, there is a natural pairing between $\mathrm{Sym}^d W$ and $\mathrm{Sym}^d W^*$. If $a = a_1 \cdots a_d \in \mathrm{Sym}^d W$, $b = b_1 \cdots b_d \in \mathrm{Sym}^d W^*$, we normalize this pairing by

$$(a,b)_d = (b,a)_d = \frac{1}{d!} \sum_{\sigma \in \mathfrak{S}_d} b_1(a_{\sigma(1)}) \cdots b_d(a_{\sigma(d)}).$$

Then if $i_1 + \cdots + i_n = d$,

$$(e_1^{i_1} \cdots e_n^{i_n}, f_1^{i_1} \cdots f_n^{i_n})_d = \frac{i_1! \cdots i_n!}{d!}.$$

Therefore, $\mathrm{Sym}^d W^*$ can be identified with the dual space of $\mathrm{Sym}^d W$.

The map
$$W \ni a \to i_d(a) = a \otimes a \otimes \cdots \otimes a \in \mathrm{Sym}^d W$$

is a polynomial map. So if $f \in \mathrm{Sym}^d W^*$, $f(a) = f(i_d(a))$ is a polynomial map from $W$ to $k$ and is homogeneous of degree $d$. We can identify $\mathrm{Sym}^d W^*$ with the space of degree $k$ forms on $W$ by this correspondence. If $a = \sum_{i=1}^n a_i e_i$,

$$i_d(a) = \sum \frac{d!}{i_1! \cdots i_n!} a_1^{i_1} \cdots a_n^{i_n} e_1^{i_1} \cdots e_n^{i_n},$$

where the sum is over all $(i_1, \cdots, i_n)$ such that $i_1 + \cdots + i_n = d$. So if $f = f_1^{i_1} \cdots f_n^{i_n}$, $f(a) = a_1^{i_1} \cdots a_n^{i_n}$. Therefore, $f$ corresponds to the monomial $a_1^{i_1} \cdots a_n^{i_n}$.

If $G \subset \mathrm{GL}(W)$ is a subgroup, $G$ acts on $W^*$ by $(gf)(v) = f(g^{-1}v)$ for $g \in G, f \in W^*$. Whenever we consider the contragredient representation, we consider this action.

## §2 Definitions of equivariant maps

Let $G, V, H$ be as in the introduction. We construct $H$–equivariant maps from $V$ or $V^{\mathrm{ss}}$ to various $H$–varieties in this section.

Let $W = k^2$ be the space of two dimensional column vectors. Let $\{e_1, e_2\}$ be the standard basis of $W$. Consider the usual action of $\mathrm{GL}(2)$ on $W$. This induces an action of $\mathrm{GL}(2)$ on $\mathrm{Sym}^d W$ and $\mathrm{Sym}^d W^*$ for any $d$. We define new actions of $\mathrm{GL}(2)$ on $\mathrm{Sym}^2 W, \mathrm{Sym}^4 W$ by $g \cdot x = (\det g)^{-1} gx, (\det g)^{-2} gx$ where $g \cdot x$ is the new action and $gx$ is the usual action. Note that scalar matrices act trivially and therefore, this defines an action of $\mathrm{PGL}(2)$ on $\mathrm{Sym}^2 W, \mathrm{Sym}^4 W$.

Let $V_1 = \mathrm{Sym}^4 W$, $V_2 = \mathrm{Sym}^2 W$, $V = \wedge^2 V_1 \otimes V_2$ and $G_1 = \mathrm{GL}(V_1) \cong \mathrm{GL}(5)$, $G_2 = \mathrm{GL}(V_2) \cong \mathrm{GL}(3)$, $G = G_1 \times G_2$. Then $G$ acts on $V$ in the usual manner and the above action of $\mathrm{PGL}(2)$ defines a homomorphism $\mathrm{PGL}(2) \to G$. By Schur's lemma, this is an imbedding. In fact, $\mathrm{Ker}(\mathrm{PGL}(2) \to G_1) = \mathrm{Ker}(\mathrm{PGL}(2) \to G_2) = \{1\}$. So we regard $\mathrm{PGL}(2)$ as a subgroup of $G$. Let $\widetilde{T} = \mathrm{Ker}(G \to \mathrm{GL}(V))$. By Schur's lemma again,

$$\widetilde{T} = \{(tI_5, t^{-2}I_3) \mid t \in \mathrm{GL}(1)\} \cong \mathrm{GL}(1).$$

If $(tI_5, t^{-2}I_3) \in \mathrm{PGL}(2)$, it acts trivially on $V_1, V_2$. So $t = 1$. Therefore, $\mathrm{PGL}(2) \cap \widetilde{T} = \{1\}$.

Let $l_0 = e_1^2, l_1 = e_1 e_2, l_2 = e_2^2$ and $m_0 = e_1^4, \cdots, m_4 = e_2^4$. Then $\{l_0, l_1, l_2\}$, $\{m_0, \cdots, m_4\}$ are bases of $V_2, V_1$ respectively.

We define a linear map $\phi_1 : V \to \wedge^4 V_1 \otimes \operatorname{Sym}^2 V_2 \cong V_1^* \otimes \operatorname{Sym}^2 V_2$ by

$$(2.1) \qquad V \ni \sum_{i=1}^N p_i \otimes q_i \to \frac{1}{2} \sum_{i,j=1}^N p_i \wedge p_j \otimes q_i q_j$$

for $p_1, \cdots, p_N \in \wedge^2 V_1$, $q_1, \cdots, q_N \in V_2$. Regarding $\operatorname{Sym}^2 V_2$ as a subspace of $V_2 \otimes V_2$, we denote the element of $V_1^* \otimes V_2 \otimes V_2$ which corresponds to $\phi_1(x)$ by $\overline{\phi}_1(x)$. Regarding $V_1^*$ as the contragredient representation of $V_1$, $V_1^* \otimes \operatorname{Sym}^2 V_1$ is a representation of $G$.

The following lemma can be proved as in [19, p. 80], and the proof is left to the reader.

**Lemma (2.2)** *For $g = (g_1, g_2) \in G$, $\phi_1(gx) = \det g_1 g \phi_1(x)$, $\overline{\phi}_1(gx) = \det g_1 g \overline{\phi}_1(x)$.*

If $x = \sum_{0 \leq i < j \leq 4} m_i \wedge m_j \otimes x_{ij}$ with $x_{ij} \in V_2$,

$$\phi_1(x) = \frac{1}{2} \sum_{0 \leq i < j \leq 4} \sum_{0 \leq k < l \leq 4} m_i \wedge m_j \wedge m_k \wedge m_l \otimes x_{ij} x_{kl}.$$

Let $m_0^*, \cdots, m_4^* \in \wedge^4 V_1$ be elements such that $m_i \wedge m_j^* = \delta_{ij} m_0 \wedge \cdots \wedge m_4$ ($\delta_{ij}$ is Kronecker's delta). Explicitly,

$$m_0 = m_1 \wedge m_2 \wedge m_3 \wedge m_4, \ m_1 = -m_0 \wedge m_2 \wedge m_3 \wedge m_4, \text{ etc.}$$

We identify $\wedge^4 V_1$ with the dual space of $V_1$ by the pairing

$$V_1 \times \wedge^4 V_1 \ni (a, b) \to a \wedge b$$

and choosing $m_0 \wedge \cdots \wedge m_4$ as the basis element of $\wedge^5 V_1$. Then $\{m_0^*, \cdots, m_4^*\}$ can be regarded as the dual basis of $\{m_0, \cdots, m_4\}$.

Let

$$(2.3) \qquad \begin{aligned} \operatorname{Pfaff}_0(x) &= x_{12} x_{34} - x_{13} x_{24} + x_{14} x_{23}, \\ \operatorname{Pfaff}_1(x) &= -(x_{02} x_{34} - x_{03} x_{24} + x_{04} x_{23}), \\ \operatorname{Pfaff}_2(x) &= x_{01} x_{34} - x_{03} x_{14} + x_{04} x_{13}, \\ \operatorname{Pfaff}_3(x) &= -(x_{01} x_{24} - x_{02} x_{14} + x_{04} x_{12}) \\ \operatorname{Pfaff}_4(x) &= x_{01} x_{23} - x_{02} x_{13} + x_{03} x_{12}. \end{aligned}$$

Then

$$(2.4) \qquad \phi_1(x) = \sum_{i=0}^4 m_i^* \otimes \operatorname{Pfaff}_i(x).$$

The quadratic polynomials $\operatorname{Pfaff}_0(x), \cdots, \operatorname{Pfaff}_4(x)$ are the Pfaffians of $4 \times 4$ main minors of $x$ if we regard $x$ as an alternating $5 \times 5$ matrix with entries in $V_1$. This idea was used in [24] for the case $G = \mathrm{GL}(5) \times \mathrm{GL}(4)$, $V = \wedge^2 k^5 \otimes k^4$ to parametrize quintic extensions of a given ground field.

Next we consider a linear map $\phi_2 : V_1^* \otimes \mathrm{Sym}^2 V_2 \to \wedge^5 \mathrm{Sym}^2 V_2 \cong \mathrm{Sym}^2 V_2^*$ defined by

$$(2.5) \qquad V_1^* \otimes \mathrm{Sym}^2 V_2 \ni \sum_{i=0}^{4} m_i^* \otimes p_i \to p_0 \wedge \cdots \wedge p_4 \in \wedge^5 \mathrm{Sym}^2 V_2 \cong \mathrm{Sym}^2 V_2^*.$$

**Definition (2.6)** $\Phi_1 = \frac{1}{3^4} \phi_2 \circ \phi_1$.

$\Phi_1$ is a map from $V$ to $\mathrm{Sym}^2 V_2^*$. The following lemma can also be proved as in [19, p. 80], and the proof is left to the reader.

**Lemma (2.7)** $\Phi_1(gx) = (\det g_1)^4 (\det g_2) g_2 \Phi_1(x)$.

This map $\Phi_1$ was also considered in [11] (using $\phi_1$ also) for a different purpose. We will show in §3 that the discriminant of $\Phi_1(x)$ is not identically zero and $V^{\mathrm{ss}}$ consists of $x$'s such that $\Phi_1(x)$ is non-degenerate.

For $x \in V$, let $\overline{\Phi}_1(x)(\alpha, \beta)$ be the symmetric bilinear form on $V_2$ associated with $\Phi_1(x)$. In other words,

$$\overline{\Phi}_1(x)(\alpha, \beta) = \frac{1}{2} (\Phi_1(x)(\alpha + \beta) - \Phi_1(x)(\alpha) - \Phi_1(x)(\beta))$$

for $\alpha, \beta \in V_2$.

We define a linear map $j_x : V_1^* \otimes V_2 \otimes V_2 \to \mathrm{Hom}(V_2, V_1^* \otimes V_2)$ by

$$(2.8) \qquad j_x(a \otimes \alpha \otimes \beta)(\gamma) = \overline{\Phi}_1(x)(\alpha, \gamma) a \otimes \beta$$

for $a \in V_1^*$, $\alpha, \beta, \gamma \in V_2$. If $f \in \mathrm{Hom}(V_2, V_1^* \otimes V_2)$, we define $gf \in \mathrm{Hom}(V_2, V_1^* \otimes V_2)$ by

$$(gf)(\alpha) = gf(g_2^{-1} \alpha)$$

where we are considering the action of $g$ on the element $f(g_2^{-1} \alpha)$.

**Lemma (2.9)** $j_{gx}(g(a \otimes \alpha \otimes \beta)) = (\det g_1)^4 (\det g_2) g j_x(a \otimes \alpha \otimes \beta)$ for all $a \in V_1^*, \alpha, \beta \in V_2$.

*Proof.* Let $\gamma \in V_2$. Then by Lemma (2.7),

$$\begin{aligned}
j_{gx}(g(a \otimes \alpha \otimes \beta))(\gamma) &= \overline{\Phi}_1(gx)(g_2 \alpha, \gamma) g_1 a \otimes g_2 \beta \\
&= (\det g_1)^4 (\det g_2)(g_2 \overline{\Phi}_1(x))(g_2 \alpha, \gamma) g_1 a \otimes g_2 \beta \\
&= (\det g_1)^4 (\det g_2) \overline{\Phi}_1(x)(\alpha, g_2^{-1} \gamma) g_1 a \otimes g_2 \beta \\
&= (\det g_1)^4 (\det g_2) g(\overline{\Phi}_1(x)(\alpha, g_2^{-1} \gamma) a \otimes \beta) \\
&= (\det g_1)^4 (\det g_2) g(j_x(a \otimes \alpha \otimes \beta)(g_2^{-1} \gamma)) \\
&= (\det g_1)^4 (\det g_2)(g(j_x(a \otimes \alpha \otimes \beta)))(\gamma).
\end{aligned}$$

This proves the lemma. □

**Definition (2.10)** $\phi_3(x) = j_x(\overline{\phi}_1(x))$.

Apparently, $\phi_3$ is a map from $V$ to $\mathrm{Hom}(V_2, V_1^* \otimes V_2)$.

**Lemma (2.11)** $\phi_3(gx) = (\det g_1)^5 (\det g_2) g \phi_3(x)$.

*Proof.*
$$\begin{aligned}
\phi_3(gx) &= j_{gx}(\overline{\phi}_1(gx)) \\
&= \det g_1 j_{gx}(g\overline{\phi}_1(x)) \\
&= (\det g_1)^5 (\det g_2) g j_x(\overline{\phi}_1(x)) \\
&= (\det g_1)^5 (\det g_2) g \phi_3(x).
\end{aligned}$$
□

By the basis $\{l_0, l_1, l_2\}$ for $V_2$, we can regard $\phi_3(x)$ as a $3 \times 3$ matrix with entries in $V_1^*$. Then the action of $g = (g_1, g_2) \in \mathrm{GL}(V_1) \times \mathrm{GL}(3)$ is obtained by considering $g_2 \phi_3(x) g_2^{-1}$ and then applying $g_1$ entry-wise. Therefore,

$$(2.12) \qquad \Phi_2(x) = \mathrm{tr}(\phi_3(x)^2) \in \mathrm{Sym}^2 V_1^*, \ F_x = \det \phi_3(x) \in \mathrm{Sym}^3 V_1^*$$

define maps $x \to \Phi_2(x), F_x$ from $V$ to $\mathrm{Sym}^2 V_1^*, \mathrm{Sym}^3 V_1^*$.

The following lemma is an easy corollary of Lemma (2.11).

**Lemma (2.13)** (1) $\Phi_2(gx) = (\det g_1)^{10} (\det g_2)^2 g_1 \Phi_2(x)$.
(2) $F_{gx}(a) = (\det g_1)^{15} (\det g_2)^3 F_x(g_1^{-1} a)$ for all $a \in V_1$.

For later purposes, we describe how to compute $\Phi_2(x), F_x$. We have already described how to compute $\phi_1(x), \Phi_1(x)$ in (2.4), (2.5). Let $\{p_0, p_1, p_2\}$ be the dual basis of $\{l_0, l_1, l_2\}$. Suppose

$$\overline{\phi}_1(x) = \sum_{i,j=0}^{2} a_{ij} \otimes l_i \otimes l_j, \ \overline{\Phi_1}(x) = \sum_{t,s=0}^{2} b_{ts} p_t \otimes p_s$$

with $a_{ij} \in V_1^*$, $b_{ts} \in k$ for all $i, j, t, s$ and $a_{ij} = a_{ji}$, $b_{ts} = b_{st}$. We denote the matrices $(a_{ij}), (b_{ts})$ also by $\overline{\phi}_1(x), \overline{\Phi}_1(x)$. Let $\gamma = \sum_{s=0}^{2} \gamma_s l_s$. Then

$$\begin{aligned}
\phi_3(x)(\gamma) &= j_x(\overline{\phi}_1(x))(\gamma) \\
&= \sum_{i,j,s=0}^{2} j_x(a_{ij} \otimes l_i \otimes l_j)(\gamma_s l_s) \\
&= \sum_{i,j,s=0}^{2} \overline{\Phi}_1(x)(l_i \otimes l_s) a_{ij} \gamma_s \otimes l_j \\
&= \sum_{i,j,s=0}^{2} a_{ij} b_{is} \gamma_s \otimes l_j \\
&= \sum_{i,j,s=0}^{2} a_{ji} b_{is} \gamma_s \otimes l_j.
\end{aligned}$$

Note that $\sum_{i=0}^{2} a_{ji}b_{is}$ is the $(j,s)$–entry of the matrix product $\overline{\phi}_1(x)\overline{\Phi}_1(x)$. So if we regard $\phi_3(x)$ as a $3 \times 3$ matrix with entries in $V_1^*$, we get the following relation

$$\phi_3(x) = \overline{\phi}_1(x)\overline{\Phi}_1(x). \tag{2.14}$$

We can express any $x \in V$ as $x = x_0 \otimes l_0 + x_1 \otimes l_1 + x_2 \otimes l_2$ where $x_0, x_1, x_2 \in \wedge^2 V_1$. Let $S_x \subset \wedge^2 V_1 \cong k^{10}$ be the subspace spanned by $x_0, x_1, x_2$. We will show in §3 that $\dim S_x = 3$ if $x \in V^{\mathrm{ss}}$. So if $x \in V^{\mathrm{ss}}$, we define

$$\Phi_3(x) \in \mathrm{Gr}(3, 10) \tag{2.15}$$

to be the element of the Grassmann determined by $S_x$. This defines a map $x \to \Phi_3(x)$ from $V^{\mathrm{ss}}$ to $\mathrm{Gr}(3, 10)$. The group $G_1$ still acts on $\mathrm{Gr}(3, 10)$ and the following lemma is obvious.

**Lemma (2.16)** $\Phi_3(gx) = g_1 \Phi_3(x)$.

**§3 Equivariant maps at $w$**

In this section, we prove that the equivariant maps we constructed in §2 are well defined and are non-trivial by evaluating them at a point $w$, which we will define in (3.10).

Let $\{l_0, l_1, l_2\}, \{m_0, \cdots, m_4\}$ be the bases of $\mathrm{Sym}^2 W, \mathrm{Sym}^4 W$ we defined in §2. Let $\{p_0, p_1, p_2\}$ be the dual basis of $\{l_0, l_1, l_2\}$. Let $\mathfrak{h}$ be the Lie algebra of $\mathrm{PGL}(2)$. Then $\mathfrak{h}$ is the Lie algebra of $\mathrm{SL}(2)$ also. We consider $V$ as a representation of $\mathfrak{h}$ also. Let $\Lambda$ be the fundamental dominant weight of $\mathfrak{h}$. We denote the irreducible representation of $\mathfrak{h}$ with highest weight $d\Lambda$ also by $d\Lambda$. Then by considering weights, $\wedge^2 V_1 \cong 6\Lambda \oplus 2\Lambda$ and

$$6\Lambda \otimes 2\Lambda \cong 8\Lambda \oplus 6\Lambda \oplus 4\Lambda, \ 2\Lambda \otimes 2\Lambda \cong 4\Lambda \oplus 2\Lambda \oplus k, \tag{3.1}$$

where $k$ is the trivial representation.

Therefore, $V$ contains the trivial representation precisely once.

Note that $V_2^* \cong \mathfrak{h}$ as an $\mathfrak{h}$–module. We identify $\mathrm{Sym}^4 W^*$ as the space of homogeneous polynomials of degree four in two variables $v = {}^t(v_1 \ v_2)$ ($v$ corresponds to $v_1 e_1 + v_2 e_2$). We identify $a = a(v) = a_0 v_1^4 + \cdots + a_4 v_2^4$ with $(a_0, \cdots, a_4)$.

Let

$$H_0 = \begin{pmatrix} 0 & 1 \\ 0 & 0 \end{pmatrix}, \ H_1 = \begin{pmatrix} 1 & 0 \\ 0 & -1 \end{pmatrix}, \ H_2 = \begin{pmatrix} 0 & 0 \\ -1 & 0 \end{pmatrix}. \tag{3.2}$$

An easy way to compute Lie algebra actions is to consider values in the ring of dual numbers $k[\epsilon]/(\epsilon^2)$.

Let

$$
(3.3) \quad A_0 = \begin{pmatrix} 0 & 1 & 0 & 0 & 0 \\ 0 & 0 & 2 & 0 & 0 \\ 0 & 0 & 0 & 3 & 0 \\ 0 & 0 & 0 & 0 & 4 \\ 0 & 0 & 0 & 0 & 0 \end{pmatrix}, \quad A'_0 = \begin{pmatrix} 0 & 1 & 0 \\ 0 & 0 & 2 \\ 0 & 0 & 0 \end{pmatrix},
$$

$$
A_1 = \begin{pmatrix} 4 & 0 & 0 & 0 & 0 \\ 0 & 2 & 0 & 0 & 0 \\ 0 & 0 & 0 & 0 & 0 \\ 0 & 0 & 0 & -2 & 0 \\ 0 & 0 & 0 & 0 & -4 \end{pmatrix}, \quad A'_1 = \begin{pmatrix} 2 & 0 & 0 \\ 0 & 0 & 0 \\ 0 & 0 & -2 \end{pmatrix},
$$

$$
A_2 = \begin{pmatrix} 0 & 0 & 0 & 0 & 0 \\ 4 & 0 & 0 & 0 & 0 \\ 0 & 3 & 0 & 0 & 0 \\ 0 & 0 & 2 & 0 & 0 \\ 0 & 0 & 0 & 1 & 0 \end{pmatrix}, \quad A'_2 = \begin{pmatrix} 0 & 0 & 0 \\ 2 & 0 & 0 \\ 0 & 1 & 0 \end{pmatrix}.
$$

Then the actions of $\mathfrak{h}$ on $V_1, V_2$ are easy to describe and with respect to the bases $\{e_0, \cdots, e_4\}, \{l_0, l_1, l_2\}$, $H_0, H_1, H_2$ map to

$$(A_0, A'_0), (A_1, A'_1), -(A_2, A'_2)$$

respectively.

Since the action of $g \in \mathrm{PGL}(2)$ on $a(v) \in \mathrm{Sym}^4 W^*$ is $a(g^{-1}v)$, the action of $H \in \mathfrak{h}$ on $a = a(v)$ is given by

$$a((1 - \epsilon H)v) = a(v) + \epsilon(Ha)(v).$$

Then easy computations show that

$$
(3.4) \quad \begin{aligned}
H_0 a &= (0, -4a_0, -3a_1, -2a_2, -a_3), \\
H_1 a &= (-4a_0, -2a_1, 0, 2a_3, 4a_4), \\
H_2 a &= (a_1, 2a_2, 3a_3, 4a_4, 0).
\end{aligned}
$$

For $a = (a_0, \cdots, a_4), b = (b_0, \cdots, b_4)$, we define

$$(3.5) \quad Q(a, b) = a_0 b_4 - \frac{1}{4} a_1 b_3 + \frac{1}{6} a_2 b_2 - \frac{1}{4} a_3 b_1 + a_4 b_0.$$

Then $Q$ is a non-degenerate symmetric bilinear form, invariant under the action of $\mathrm{PGL}(2)$. This implies $Q(Ha, b)$ is an alternating form for any $H \in \mathfrak{h}$. We regard this alternating form as an element of $\wedge^2(\mathrm{Sym}^4 W^*)^* \cong \wedge^2 \mathrm{Sym}^4 W$.

**Lemma (3.6)** *The map $H \to f_H(a, b) = Q(Ha, b)$ is an $\mathfrak{h}$–homomorphism from $\mathfrak{h}$ to $\wedge^2 \mathrm{Sym}^4 W$.*

*Proof.* Note that the action of $H \in \mathfrak{h}$ on an element $f(a,b)$ in $\wedge^2(\mathrm{Sym}^4 W^*)^*$ is given by $(Hf)(a,b) = -f(Ha,b) - f(a,Hb)$. So if $H, H' \in \mathfrak{h}$,

$$\begin{aligned}
(H'f_H)(a,b) &= -f_H(H'a,b) - f_H(a, H'b) \\
&= -Q(HH'a, b) - Q(Ha, H'b) \\
&= -Q(HH'a, b) + Q(H'Ha, b) \\
&= Q([H', H]a, b) \\
&= f_{[H', H]}(a, b).
\end{aligned}$$

□

This defines an $\mathfrak{h}$–homomorphism $\mathfrak{h} \to \wedge^2 V_1$. Regarding this homomorphism as an element of $\wedge^2 V_1 \otimes \mathfrak{h}^* \cong \wedge^2 V_1 \otimes V_2 = V$, we get a fixed point of $V$ under the action of $\mathrm{PGL}(2)$. We compute this element explicitly.

Note that the linear map defined by

(3.7) $$H_0 \to \frac{1}{2}v_2^2, \ H_1 \to v_1 v_2, \ H_2 \to \frac{1}{2}v_1^2$$

is an $\mathfrak{h}$–homomorphism.

By (3.4),

(3.8) $$\begin{aligned}
Q(H_0 a, b) &= a_0 b_3 - \frac{1}{2} a_1 b_2 + \frac{1}{2} a_2 b_1 - a_3 b_0, \\
Q(H_1 a, b) &= -4 a_0 b_4 + \frac{1}{2} a_1 b_3 - \frac{1}{2} a_3 b_1 + 4 a_4 b_0, \\
Q(H_2 a, b) &= a_1 b_4 - \frac{1}{2} a_2 b_3 + \frac{1}{2} a_3 b_2 - a_4 b_1.
\end{aligned}$$

Note that $\{e_1^2, 2e_1 e_2, e_2^2\}$ is the dual basis of $\{v_1^2, v_1 v_2, v_2^2\}$, and $(m_0, v_1^4)_4 = 1$, $(m_1, v_1^3 v_2)_4 = \frac{1}{4}$, etc. We identify $\wedge^2 V_1$ with the space of alternating bilinear forms on $V_1^*$ by assuming

$$m \wedge m'(a,b) = (m, a)_1 (m', b)_1 - (m, b)_1 (m', a)_1$$

for $m, m' \in V_1$, $a, b \in V_1^*$. So by the correspondence $H \to f_H$,

(3.9) $$\begin{aligned}
H_0 &\to 4 m_0 \wedge m_3 - 12 m_1 \wedge m_2, \\
H_1 &\to -4 m_0 \wedge m_4 + 8 m_1 \wedge m_3, \\
H_2 &\to 4 m_1 \wedge m_4 - 12 m_2 \wedge m_3.
\end{aligned}$$

Since $\{H_0, H_1, H_2\}$ corresponds to $\{\frac{1}{2}v_2^2, v_1 v_2, \frac{1}{2}v_1^2\}$ and $\{2l_2, 2l_1, 2l_0\}$ is its dual basis, this correspondence can be regarded as the element $8w$ where

(3.10) $$\begin{aligned}
w = &(m_0 \wedge m_3 - 3 m_1 \wedge m_2) \otimes l_2 + (-m_0 \wedge m_4 + 2 m_1 \wedge m_3) \otimes l_1 \\
&+ (m_1 \wedge m_4 - 3 m_2 \wedge m_3) \otimes l_0.
\end{aligned}$$

These considerations show the following proposition.

**Proposition (3.11)** *The element $w \in V$ is fixed by $\mathrm{PGL}(2)$.*

In [19, p. 95], instead of $w$, the element

$$w' = (m_0 \wedge m_1 + m_2 \wedge m_3, m_1 \wedge m_2 + m_3 \wedge m_4, m_0 \wedge m_2 + m_1 \wedge m_4)$$

was considered (we shifted the indices in [19] because we are using indices $0, \cdots, 4$). However, by replacing $m_0, \cdots, m_4$ in $w$ by $m_4, m_2, m_0, m_1, m_3$ respectively, and multiplying scalars to basis elements of $\wedge^2 V_1$, we get the above element $w'$. Therefore, we are considering essentially the same element as in [19].

In [19, p. 96], the Lie algebra of $G^0_{w'}/\widetilde{T}$ ($G^0_{w'}$ is the connected component of the stabilizer) is computed and is isomorphic to the Lie algebra of $\mathrm{PGL}(2)$. Therefore, this is the case for $w$ also. Since we are assuming $\mathrm{ch}\, k = 0$, this implies $G^0_w = \mathrm{PGL}(2) \times \widetilde{T}$ if $k$ is algebraically closed (see [7]). For arbitrary $k$, we still have the inclusion $\mathrm{PGL}(2) \times \widetilde{T} \subset G^0_w$. Since this is an isomorphism over $\bar{k}$, we get the following proposition

**Proposition (3.12)** $G^0_w = \mathrm{PGL}(2) \times \widetilde{T}$.

By the basis $\{e_0, \cdots, e_4\}$, we regard $V$ as the space of $5 \times 5$ alternating matrices with entries in $V_2$. Then

$$w = \begin{pmatrix} 0 & 0 & 0 & l_2 & -l_1 \\ 0 & 0 & -3l_2 & 2l_1 & l_0 \\ 0 & 3l_2 & 0 & -3l_0 & 0 \\ -l_2 & -2l_1 & 3l_0 & 0 & 0 \\ l_1 & -l_0 & 0 & 0 & 0 \end{pmatrix}.$$

By the definition (2.3),

(3.13) $\mathrm{Pfaff}_0(w) = -3l_0^2$, $\mathrm{Pfaff}_1(w) = -3l_0 l_1$, $\mathrm{Pfaff}_2(w) = -2l_1^2 - l_0 l_2$,
$\mathrm{Pfaff}_3(w) = -3l_1 l_2$, $\mathrm{Pfaff}_4(w) = -3l_2^2$.

Note that we are regarding them as elements of $\mathrm{Sym}^2 V_2$ and not $\mathrm{Sym}^4 W$.

By the basis $\{l_0, l_1, l_2\}$, we regard $\overline{\phi}_1(w)$ as a $3 \times 3$ matrix with entries in $V_1^*$ as in §2. Then

(3.14) $$\overline{\phi}_1(w) = - \begin{pmatrix} 3m_0^* & \frac{3}{2}m_1^* & \frac{1}{2}m_2^* \\ \frac{3}{2}m_1^* & 2m_2^* & \frac{3}{2}m_3^* \\ \frac{1}{2}m_2^* & \frac{3}{2}m_3^* & 3m_4^* \end{pmatrix}.$$

Let

(3.15) $n_0 = l_0^2$, $n_1 = l_1^2$, $n_2 = l_2^2$, $n_3 = l_0 l_1$, $n_4 = l_1 l_2$, $n_5 = l_0 l_2 \in \mathrm{Sym}^2 V_2$.

Then $\{n_0, \cdots, n_5\}$ is a basis of $\mathrm{Sym}^2 V_2$. Let $n_0^*, \cdots, n_5^*$ be elements of $\wedge^5 \mathrm{Sym}^2 V_2 \cong \mathrm{Sym}^2 V_2^*$ such that $n_i \wedge n_j^* = \delta_{ij} n_0 \wedge \cdots \wedge n_5$. Then

$$\mathrm{Pfaff}_0(w) \wedge \mathrm{Pfaff}_1(w) \wedge \mathrm{Pfaff}_2(w) \wedge \mathrm{Pfaff}_3(w) \wedge \mathrm{Pfaff}_4(w)$$
$$= (-3n_0) \wedge (-3n_3) \wedge (-2n_1 - n_5) \wedge (-3n_4) \wedge (-3n_2)$$
$$= -3^4 n_0 \wedge n_3 \wedge (2n_1 + n_5) \wedge n_4 \wedge n_2$$
$$= -3^4 (2n_0 \wedge n_3 \wedge n_1 \wedge n_4 \wedge n_2 + n_0 \wedge n_3 \wedge n_5 \wedge n_4 \wedge n_2)$$
$$= 3^4 (n_1^* - 2n_5^*).$$

Therefore,

(3.15) $$\Phi_1(w) = n_1^* - 2n_5^*.$$

We identify $\mathrm{Sym}^2 V_2^*$ with the dual space of $\mathrm{Sym}^2 V_2$. Then with respect to the basis $\{p_1 p_1, p_2\}$, $n_1^*, n_5^*$ correspond to $p_1^2, 2p_0 p_2$. Therefore, $\Phi_1(w) = p_1^2 - 4p_0 p_2$. We regard $\overline{\Phi}_1(w)$ as a $3 \times 3$ matrix as in §2. Then

(3.16) $$\overline{\Phi}_1(w) = \begin{pmatrix} 0 & 0 & -2 \\ 0 & 1 & 0 \\ -2 & 0 & 0 \end{pmatrix}.$$

Therefore,

(3.17) $$\phi_3(w) = \overline{\phi}_1(w)\overline{\Phi}_1(w) = \begin{pmatrix} m_2^* & -\frac{3}{2}m_1^* & 6m_0^* \\ 3m_3^* & -2m_2^* & 3m_1^* \\ 6m_4^* & -\frac{3}{2}m_3^* & m_2^* \end{pmatrix}.$$

So, we get the following proposition easily.

**Proposition (3.18)** *Let $a = a_0 m_0 + \cdots + a_4 m_4$. Then*
(1) $\Phi_2(w)(a) = (\mathrm{tr}(\phi_3(w)^2))(a) = 72(a_0 a_4 - \frac{1}{4} a_1 a_3 + \frac{1}{12} a_2^2),$
(2) $F_w(a) = (\det \phi_3(w))(a) = 72 a_0 a_2 a_4 + 9 a_1 a_2 a_3 - 2a_2^3 - 27 a_0 a_3^2 - 27 a_1^2 a_4.$

By these considerations, $\Phi_1, \Phi_2$ are non-trivial maps. By (3.16), the discriminant $\Delta(x)$ of $\Phi_1(x)$ is a non-zero polynomial. By Lemma (2.7), $\Delta(x)$ is a non-constant polynomial. Therefore, it is a relative invariant polynomial. So we reproved that $w \in V_k^{\mathrm{ss}}$. Since our case is known to be a regular prehomogeneous vector space, $V_k^{\mathrm{ss}}$ is a single $G_k$–orbit if $k$ is algebraically closed. Therefore, $V^{\mathrm{ss}}$ consists of $x$'s such that $\Phi_1(x)$ is non-degenerate.

Since $\Phi_2(w)$ is non-degenerate, $\Phi_2(x)$ is non-degenerate for all $x \in V^{\mathrm{ss}}$. Apparently, three components of $w$ are linearly independent. Therefore, $\Phi_3(x) \in \mathrm{Gr}(3, 10)$ is well defined for all $x \in V^{\mathrm{ss}}$ also.

## §4 The orbit space $G_k \setminus V_k^{\mathrm{ss}}$

In this section, we prove that $G_k \setminus V_k^{\mathrm{ss}}$ corresponds bijectively with $\mathrm{GL}(1)_k \times \mathrm{GL}(3)_k$–equivalence classes of ternary quadratic forms over $k$.

We first recall the relation between the orbit space $G_k \setminus V_k^{\mathrm{ss}}$ and the Galois cohomology set.

For any algebraic group $G$ over $k$, let $\mathrm{H}^1(k, G)$ be the first Galois cohomology set. We choose the definition so that trivial classes are those of the form $\{g^{-1}g^\sigma\}_{\sigma \in \mathrm{Gal}(\bar k/k)}$ ($g \in G_{\bar k}$) and the cocycle condition is $h_{\sigma\tau} = h_\tau h_\sigma^\tau$ for a continuous map $\{h_\sigma\}_{\sigma \in \mathrm{Gal}(\bar k/k)}$ from $\mathrm{Gal}(\bar k/k)$ to $G_{\bar k}$.

Let $(G, V)$ be an arbitrary regular prehomogeneous vector space, and $w \in V_k^{\mathrm{ss}}$. Then for any $x \in V_k^{\mathrm{ss}}$, there exists $g_x \in G_{\bar k}$ such that $x = g_x w$. Then $c_x = \{g_x^{-1}g_x^\sigma\}_{\sigma \in \mathrm{Gal}(\bar k/k)}$ determines a cohomology class in $\mathrm{H}^1(k, G_w)$ and does not depend on the choice of $g_x$. The following theorem is due to Igusa [8].

**Theorem (4.1) (Igusa)** *The correspondence*

$$G_k \setminus V_x^{\mathrm{ss}} \ni x \to c_x \in \mathrm{Ker}(\mathrm{H}^1(k, G_w) \to \mathrm{H}^1(k, G))$$

*is bijective.*

Note that $\mathrm{Ker}(\mathrm{H}^1(k, G_w) \to \mathrm{H}^1(k, G))$ is the set of elements $c \in \mathrm{H}^1(k, G_w)$ which map to the trivial class in $\mathrm{H}^1(k, G)$. In our case, $\mathrm{H}^1(k, G)$ is trivial. Therefore, $G_k \setminus V_k^{\mathrm{ss}} \cong \mathrm{H}^1(k, G_w)$.

We recall the correspondence between $\mathrm{GL}(1)_k \times \mathrm{GL}(3)_k$–equivalence classes of ternary quadratic forms and quarternion algebras. Let $V_3 = \mathrm{Sym}^2 V_2^*$. Then $\mathrm{GL}(V_2) \cong \mathrm{GL}(3)$ acts on $V_3$ in the usual manner. We let $\mathrm{GL}(1)$ act on $V_3$ by the usual multiplication. Then $\mathrm{GL}(1) \times \mathrm{GL}(3)$ acts on $V_3$. Let $\{l_0, l_1, l_2\}$ and $\{p_0, p_1, p_2\}$ be as before. Let $\overline{w} = p_1^2 - 4p_0 p_2$, and $Q$ the corresponding quadratic form. It is well known (and is easy to verify) that the stabilizer of $\overline{w}$ is isomorphic to $\mathrm{SO}(Q) \times \mathrm{GL}(1)$ and $\mathrm{SO}(Q) \cong \mathrm{PGL}(2)$.

Therefore, $(\mathrm{GL}(1)_k \times \mathrm{GL}(3)_k) \setminus V_{3k}^{\mathrm{ss}}$ corresponds bijectively with $\mathrm{H}^1(k, \mathrm{PGL}(2))$. Since $\mathrm{PGL}(2)$ is isomorphic to the automorphism group of the associative algebra $\mathrm{M}(2,2)$, $\mathrm{H}^1(k, \mathrm{PGL}(2))$ corresponds bijectively with isomorphism classes of quarternion algebras. Given a ternary quadratic form, the corresponding quarternion algebra is the Clifford algebra associated with the quadratic form.

Now we go back to our situation. Let $G, H, V$ be as before. We consider the element $w \in V_k^{\mathrm{ss}}$ which we defined in (3.10). We pointed out in (3.12) that $G_w^0 \cong \mathrm{PGL}(2) \times \mathrm{GL}(1)$.

**Proposition (4.2)** *The group $G_w$ is connected.*

*Proof.* We may assume that $k$ is algebraically closed. Suppose $g \in G_k/\widetilde{T}_k$. The connected component of the stabilizer of $w$ in $G/\widetilde{T}$ is isomorphic to $\mathrm{PGL}(2)$. The conjugation by $g$ induces an automorphism of $\mathrm{PGL}(2)$. Since there is no outer automorphism of $\mathrm{PGL}(2)$, by changing $g$ if necessary, we may assume that $g$ commutes with elements of $\mathrm{PGL}(2)$. Since $V_1, V_2$ are irreducible representations, by Schur's lemma, $g$ is represented by an element of the form $(t_1 I_5, t_2 I_3)$. This element fixes $w$ if and only if $t_1^2 t_2 = 1$. So $g = 1$ (in $G/\widetilde{T}$). $\square$

**Proposition (4.3)** *(1) The map $\Phi_1 : V \to V_3 = \mathrm{Sym}^2 V_1^*$ induces a bijection $G_k \setminus V_k^{\mathrm{ss}} \cong (\mathrm{GL}(1)_k \times \mathrm{GL}(3)_k) \setminus V_{3k}^{\mathrm{ss}}$.*
*(2) If $x \in V_k^{\mathrm{ss}}$, the projections of $H_x$ to $G_1, G_2$ induce isomorphisms to the images. In particular, $H_x \cong \mathrm{SO}(\Phi_1(x))$.*

*Proof.* Let $c \in \mathrm{H}^1(k, G_w)$. Then $c$ becomes trivial in $\mathrm{H}^1(k, H)$ also. Let $g = (g_1, g_2) \in H_{\bar{k}}$ be the element such that $c$ is represented by $\{g^{-1}g^\sigma\}_{\sigma \in \mathrm{Gal}(\bar{k}/k)}$. Then the orbit in $V_k^{\mathrm{ss}}$ corresponding to $c$ is $gw$. By Lemma (2.7), $\Phi_1(gw) = g_2 \Phi_1(w)$.

Since $\mathrm{H}^1(k, G_w) \cong \mathrm{H}^1(k, \mathrm{PGL}(2))$ and the projection of $\mathrm{PGL}(2)$ to $G_2$ is an isomorphism to its image,

$$\mathrm{H}^1(k, G_w) \ni \{g^{-1}g^\sigma\}_{\sigma \in \mathrm{Gal}(\bar{k}/k)} \to \{g_2^{-1}g_2^\sigma\}_{\sigma \in \mathrm{Gal}(\bar{k}/k)} \in \mathrm{H}^1(k, \mathrm{PGL}(2))$$

is a bijection. Note that we are considering $\mathrm{PGL}(2) \subset G$ for the first element and $\mathrm{PGL}(2) \subset G_2$ for the second element. Since $(\mathrm{GL}(1)_k \times \mathrm{GL}(3)_k) \setminus V_{3k}^{\mathrm{ss}} \cong \mathrm{H}^1(k, \mathrm{PGL}(2))$, this proves (1).

Note that $(\mathrm{PGL}(2) \times \widetilde{T}) \cap H = \mathrm{PGL}(2)$. So $H_w \cong \mathrm{PGL}(2)$. We already pointed out that statement (2) holds for $w$ in §2. Let $x \in V_k^{\mathrm{ss}}$. By Lemma (2.7), the projection of $H_x$ to $G_2$ is contained in $\mathrm{SO}(\Phi_1(x))$. So it is enough to prove (2) when $k$ is algebraically closed. But then $x$ is in the orbit of $w$ and (2) follows easily. □

**Remark (4.4)** The map $\Phi_1$ induces a map $G_k \setminus V_k^{\mathrm{ss}} \to \mathrm{GL}(3)_k \setminus V_{3k}^{\mathrm{ss}}$ also, but this may not be surjective. This may be regarded as the section $(\mathrm{GL}(1)_k \times \mathrm{GL}(3)_k) \setminus V_{3k}^{\mathrm{ss}} \to \mathrm{GL}(3)_k \setminus V_{3k}^{\mathrm{ss}}$ defined by $x \to (\det x)^{-1}x$.

### §5 Intermediate groups

Let $x \in V_{\mathbb{R}}^{\mathrm{ss}}$. By Proposition (4.3), $H_{x\mathbb{R}}$ is connected in classical topology. So $H_{x\mathbb{R}+}^0 = H_{x\mathbb{R}}$. If $\Phi_1(x)$ is definite, $H_{x\mathbb{R}}$ is compact by Proposition (4.3) also. Then $H_{x\mathbb{R}} H_{\mathbb{Z}} \subset H_{\mathbb{R}}/H_{\mathbb{Z}}$ is a compact set. Therefore, an analogue of the Oppenheim conjecture is not applicable to such points. The set of real indefinite non-degenerate ternary quadratic forms is a single $\mathrm{GL}(1)_{\mathbb{R}} \times \mathrm{GL}(3)_{\mathbb{R}}$–orbit. Therefore, we only consider $x \in G_{\mathbb{R}} w$.

We determine all the closed connected subgroups between $H_{x\mathbb{R}+}^0$ and $H_{\mathbb{R}}$ for all $x \in G_{\mathbb{R}} w$ for the rest of this section. This reduces to the consideration of Lie algebras. We consider an arbitrary ground field $k$ of characteristic zero and specialize to $k = \mathbb{R}$ in (5.10).

We first describe possible candidates for such subgroups. By Lemmas (2.7), (2.13), (2.16), $\Phi_1, \Phi_2, \Phi_3$ are $H$–equivariant maps. As we pointed out at the end of §3, $\Phi_1(x) \in \mathrm{Sym}^2 V_1^*$, $\Phi_2(x) \in \mathrm{Sym}^2 V_2^*$ are non-degenerate for $x \in G_k w$. So let $\mathrm{SO}(\Phi_1(x)), \mathrm{SO}(\Phi_2(x))$ be the corresponding special orthogonal groups.

In the following definition, $x \in G_k w$.

**Definition (5.1)** (1) $H_{x1} \subset \mathrm{GL}(V_1), H_{x2} \subset \mathrm{GL}(V_2)$ *are the images of the projections of $H_x$ to $G_1, G_2$ respectively.*
(2) $H_{x3} = \mathrm{SO}(\Phi_2(x)) \subset \mathrm{GL}(V_1)$.

Note that both $H_{x1}, H_{x2}$ are isomorphic to $H_x$, and $H_x \cong \mathrm{PGL}(2)$.

Let $\mathfrak{h}$ be the Lie algebra of $\mathrm{PGL}(2)$ as before. Let $\mathfrak{h}_1 = \mathrm{sl}(5), \mathfrak{h}_2 = \mathrm{sl}(3)$ (Lie algebras of $\mathrm{SL}(5), \mathrm{SL}(3)$). If $\mathfrak{f}$ is a Lie algebra between $\mathfrak{h}$ and $\mathfrak{h}_1 \times \mathfrak{h}_2$, it it is an $\mathfrak{h}$–module. So we first decompose $\mathfrak{h}_1, \mathfrak{h}_2$ to direct sums of irreducible $\mathfrak{h}$–modules.

Let

(5.2) $$B = B(b_0, \cdots, b_4) = \begin{pmatrix} 2b_2 & -3b_1 & b_0 & 0 & 0 \\ 12b_3 & -b_2 & -2b_1 & 3b_0 & 0 \\ 6b_4 & 3b_3 & -2b_2 & 3b_1 & 6b_0 \\ 0 & 3b_4 & -2b_3 & -b_2 & 12b_1 \\ 0 & 0 & b_4 & -3b_3 & 2b_2 \end{pmatrix},$$

$$C = C(c_0, \cdots, c_6) = \begin{pmatrix} c_3 & 3c_2 & -c_1 & c_0 & 0 \\ 12c_4 & -2c_3 & -4c_2 & 0 & 4c_0 \\ 6c_5 & -6c_4 & 0 & -6c_2 & 6c_1 \\ 4c_6 & 0 & -4c_4 & 2c_3 & 12c_2 \\ 0 & c_6 & -c_5 & 3c_4 & -c_3 \end{pmatrix},$$

$$D = D(d_0, \cdots, d_8) = \begin{pmatrix} d_4 & -d_3 & d_2 & -d_1 & d_0 \\ 4d_5 & -4d_4 & 4d_3 & -4d_2 & 4d_1 \\ 6d_6 & -6d_5 & 6d_4 & -6d_3 & 6d_2 \\ 4d_7 & -4d_6 & 4d_5 & -4d_4 & 4d_3 \\ d_8 & -d_7 & d_6 & -d_5 & d_4 \end{pmatrix},$$

$$B'(b'_0, \cdots, b'_4) = \begin{pmatrix} b'_2 & -b'_1 & b'_0 \\ 2b'_3 & -2b'_2 & 2b'_1 \\ b'_4 & -b'_3 & b'_2 \end{pmatrix},$$

where $b_0 \cdots \in k$.

We define

(5.3)  $U_2 = \{B(b_0, \cdots, b_4) \mid b_0, \cdots, b_4 \in k\},$
$U_3 = \{C(c_0, \cdots, c_6) \mid c_0, \cdots, c_6 \in k\},$
$U_4 = \{D(d_0, \cdots, d_8) \mid d_0, \cdots, d_8 \in k\},$
$V_2 = \{B'(b'_0, \cdots, b'_4) \mid b'_0, \cdots, b'_4 \in k\}.$

Let $U_1, V_1$ be the images of $\mathfrak{h}$ in $\mathfrak{h}_1, \mathfrak{h}_2$. $U_1, V_1$ are clearly, sub $\mathfrak{h}$–modules.

**Lemma (5.4)** *The subspaces $U_2, U_3, U_4, V_2$ are irreducible sub $\mathfrak{h}$–modules with highest weights $4\Lambda, 6\Lambda, 8\Lambda, 4\Lambda$ respectively.*

*Proof.* By straightforward computations,

(5.5) $[A_0, B(1, 0, \cdots, 0)] = [A_0, C(1, 0, \cdots, 0)] = [A_0, D(1, 0, \cdots, 0)] = 0,$
$[A'_0, B'(1, 0, \cdots, 0)] = 0,$
$[A_1, B(1, 0, \cdots, 0)] = 4B(1, 0, \cdots, 0),$
$[A_1, C(1, 0, \cdots, 0)] = 6C(1, 0, \cdots, 0),$
$[A_1, D(1, 0, \cdots, 0)] = 8D(1, 0, \cdots, 0),$
$[A'_1, B'(1, 0, \cdots, 0)] = 4B'(1, 0, \cdots, 0).$

Also,

(5.6)  $[A_2, B(b_0, \cdots, b_4)] = B(0, b_0, 6b_1, b_2, 4b_3),$
$[A_2, C(c_0, \cdots, c_6)] = C(0, 2c_0, c_1, -12c_2, c_3, 10c_4, 3c_5),$
$[A_2, D(d_0, \cdots, d_8)] = D(0, d_0, 2d_1, 3d_2, \cdots, 8d_7),$
$[A'_2, B'(b'_0, \cdots, b'_4)] = B'(0, b'_0, 2b'_1, 3b'_2, 4b'_3).$

The author used MAPLE [1] to find $U_2, U_3, U_4$ but computed (5.5), (5.6) manually and checked (5.5), (5.6) by MAPLE again.

By (5.6), $U_2$ is spanned by elements of the form $\mathrm{ad}(A_2)^i B(1,0,0,0,0)$ ($\mathrm{ad}(*)$ is the adjoint representation). Since

$$\mathrm{ad}(A_0)B(1,0,0,0,0) = 0, \ \mathrm{ad}(A_1)B(1,0,0,0,0) = 4B(1,0,0,0,0),$$

$U_2$ is an irreducible sub $\mathfrak{h}$–module with highest weight $4\Lambda$.

Other cases are similar. $\square$

**Proposition (5.7)** (1) $[U_2, U_2] = U_1 \oplus U_3$.

(2) $[U_2, U_3] = U_2 \oplus U_4$.

(3) $[U_2, U_4] = U_3$.

(4) $[U_3, U_3] = U_1 \oplus U_3$.

(5) $[U_3, U_4] = U_2 \oplus U_4$.

(6) $[U_4, U_4] = U_1 \oplus U_3$.

(7) $[V_2, V_2] = V_1$.

*Proof.* We first consider (1). Since $U_2$ is irreducible, for any non-zero element $X \in U_2$, $U_2$ is generated by $X$ as an $\mathfrak{h}$–module. So $[U_2, U_2]$ is generated by $[X, U_2]$ as an $\mathfrak{h}$–module also.

By straightforward computations,

(5.8) $$[B(0,0,1,0,0), B(b_0, \cdots, b_4)]$$
$$= -\frac{21b_1}{5}A_0 - \frac{21b_3}{5}A_2$$
$$+ C(0, -4b_0, -\frac{8b_1}{5}, 0, -\frac{8b_3}{5}, -4b_4, 0).$$

We chose $B(0,0,1,0,0)$ because it is diagonal.

By (5.8), $[U_2, U_2] \subset U_1 \oplus U_3$. By choosing $b_1 = b_3 = 0$ in (5.8), $[U_2, U_2]$ contains a non-zero element of $U_3$. This implies $[U_2, U_2]$ contains $U_3$. By choosing $b_1 \neq 0$ in (5.8), $[U_2, U_2]$ contains an element of the form $X + X'$ where $X \in U_1$ is non-zero and $X' \in U_3$. So $X \in [U_2, U_2]$. This implies $[U_2, U_2]$ contains $U_1$ also. This proves (1).

Other cases follow from the following relations and by similar arguments. We found these relations manually. However, it can be checked by a routine program in MAPLE (which we did).

$$[B(0,0,1,0,0), C(c_0, \cdots, c_6)]$$
$$= B(-\frac{16c_1}{7}, -\frac{16c_2}{7}, 0, -\frac{16c_4}{7}, -\frac{16c_5}{7})$$
$$\quad + D(0, -3c_0, -\frac{12c_1}{7}, -\frac{15c_2}{7}, 0, -\frac{15c_4}{7}, -\frac{12c_5}{7}, -3c_6, 0),$$
$$[B(0,0,1,0,0), D(d_0, \cdots, d_8)]$$
$$= C(-3d_1, -4d_2, -d_3, 0, -d_5, -4d_6, -3d_7),$$
$$[C(0,0,0,1,0,0,0), C(c_0, \cdots, c_6)],$$
$$= 6c_2 A_0 - 6c_4 A_2 + C(-c_0, c_1, c_2, 0, -c_4, -c_5, c_6),$$
$$[C(0,0,0,1,0,0,0), D(d_0, \cdots, d_8)]$$
$$= B(\frac{20d_2}{7}, \frac{10d_3}{7}, 0, -\frac{10d_5}{7}, -\frac{20d_6}{7})$$
$$\quad + D(2d_0, -d_1, -\frac{13d_2}{7}, -\frac{9d_3}{7}, 0, \frac{9d_5}{7}, \frac{13d_6}{7}, d_7, -2d_8),$$
$$[D(0,0,0,0,1,0,0,0,0), D(d_0, \cdots, d_8)]$$
$$= -14d_3 A_0 - 14d_5 A_2 + C(-5d_1, 5d_2, 3d_3, 0, 3d_5, 5d_6, -5d_7),$$
$$[B'(0,0,1,0,0), B'(b'_0, \cdots, b'_4)]$$
$$= -3b_1 A'_0 - 3b_3 A'_2.$$

□

Note that the Lie algebras of $H_{w1}, H_{w2}$ are isomorphic to sl(2). We denote the Lie algebra of $H_{w3}$ by so(5) (more precisely so(3,2)). Since $\dim \text{so}(5) = 10$, $\text{so}(5) = U_1 \oplus U_3$ by counting the dimension.

**Proposition (5.9)** *If $\mathfrak{h} \subset \mathfrak{f} \subset \mathfrak{h}_1 \times \mathfrak{h}_2$ is a Lie subalgebra, $\mathfrak{h}$ is one of the following subalgebras.*

$$\mathfrak{h}, \text{ sl}(2) \times \text{sl}(2), \text{ sl}(2) \times \text{sl}(3), \text{ so}(5) \times \text{sl}(2),$$
$$\text{so}(5) \times \text{sl}(3), \text{ sl}(5) \times \text{sl}(2), \text{ sl}(5) \times \text{sl}(3).$$

*Proof.* Let $\mathfrak{f}$ be as above. Note that sl(3) does not contain any $\mathfrak{h}$–module which is isomorphic to $U_3$ or $U_4$. Suppose $\mathfrak{f} \supset U_4$. Then $\mathfrak{f} \supset U_1 \oplus U_3$ by Lemma (5.7)(1). So $\mathfrak{f} \supset U_2$ by Lemma (5.7)(5). Since $\mathfrak{f} \supset U_1$, $\mathfrak{f} \supset U_1 \oplus V_1$. Therefore, $\mathfrak{f} \supset \text{sl}(5) \times \text{sl}(2)$. So $\mathfrak{f} = \text{sl}(5) \times \text{sl}(2)$ or $\text{sl}(5) \times \text{sl}(3)$.

Suppose the projection of $\mathfrak{f}$ to the first factor contains $U_2$. Then there exists an $\mathfrak{h}$–homomorphism $\alpha : U_2 \to V_2$ such that $(x, \alpha(x)) \in \mathfrak{f}$ for all $x \in U_2$. By Lemma (5.7)(1), the projection of $\mathfrak{f}$ to the first factor contains $U_3$. Since $U_3$ is not equivalent to any other factor, $\mathfrak{f}$ contains $U_3$. By Lemma (5.7)(4), $\mathfrak{f} \supset U_1$. Since $\mathfrak{f} \supset \mathfrak{h}$, $\mathfrak{f} \supset U_1 \oplus V_1$. If $x \in U_2, y \in U_1$, $(y, 0) \in \mathfrak{f}$. So $[(y,0), (x, \alpha(x))] = ([y,x], 0) \in \mathfrak{f}$. Since $[U_1, U_2] = U_2$, $\mathfrak{f} \supset U_2$. By Lemma (5.7)(2), $\mathfrak{f} \supset U_4$ and it reduces to the previous case.

Suppose $\mathfrak{f}$ does not contain $U_4$ and the projection to the first factor does not contain $U_2$. Suppose $\mathfrak{f} \supset U_3$. By Lemma (5.7)(4), $\mathfrak{f} \supset U_1$. Therefore, $\mathfrak{f}$ has so(5) as the first factor. This implies $\mathfrak{f} = \text{so}(5) \times \text{sl}(2)$ or $\text{so}(5) \times \text{sl}(3)$.

Suppose the projection of $\mathfrak{f}$ to the first factor is $U_1$. If $\mathfrak{f} \supset V_2$, $\mathfrak{f} \supset V_1$ also. Therefore, $\mathfrak{f} = \mathrm{sl}(2) \times \mathrm{sl}(3)$. Otherwise the projection of $\mathfrak{f}$ to both factors are $\mathrm{sl}(2)$. Since there is no sub $\mathfrak{h}$–module between $\mathfrak{h}$ and $U_1 \times V_1$, $\mathfrak{f}$ is $\mathfrak{h}$ or $\mathrm{sl}(2) \times \mathrm{sl}(2)$. $\square$

Now we specialize to the field $k = \mathbb{R}$.

**Proposition (5.10)** Let $x \in G_\mathbb{R} w$ and $H_{x\mathbb{R}} \subset F \subset H_\mathbb{R}$ be a closed connected subgroup. Then $F$ is one of the following subgroups.

$$H_{x\mathbb{R}},\ H_{x1\mathbb{R}} \times H_{x2\mathbb{R}},\ H_{x1\mathbb{R}} \times \mathrm{SL}(3)_\mathbb{R},$$
$$H_{x3\mathbb{R}} \times H_{x2\mathbb{R}},\ H_{x3\mathbb{R}} \times \mathrm{SL}(3)_\mathbb{R},\ \mathrm{SL}(5)_\mathbb{R} \times H_{x2\mathbb{R}},\ \mathrm{SL}(5)_\mathbb{R} \times \mathrm{SL}(3)_\mathbb{R}.$$

*Proof.* If $x = gw$ for $g \in G_\mathbb{R}$, $H_{x\mathbb{R}} = gH_{w\mathbb{R}}g^{-1}$, etc. So we may assume that $x = w$. Then this proposition follows from the previous proposition. $\square$

### §6 An analogue of the Oppenheim conjecture

In this section, we prove an analogue of the Oppenheim conjecture.
In the following lemma, $x \in V_\mathbb{C}^{\mathrm{ss}}$. We define $H_{x1\mathbb{C}}$, etc. as in Definition (5.1).

**Lemma (6.1)** (1) If $y \in V_\mathbb{C}$ is fixed by $H_{x\mathbb{C}}$, $y$ is a scalar multiple of $x$.
(2) If $y \in \mathrm{Gr}(3, 10)_\mathbb{C}$ is fixed by $H_{x1\mathbb{C}}$, $y = \Phi_3(x)$.
(3) If $y \in \mathrm{Sym}^2 V_2^*$ is fixed by $H_{x2\mathbb{C}}$, $y$ is a scalar multiple of $\Phi_1(x)$.
(4) If $y \in \mathrm{Sym}^2 V_1^*$ is fixed by $H_{x3\mathbb{C}}$, $y$ is a scalar multiple of $\Phi_2(x)$.

*Proof.* Consider (1). Let $x = gw$ with $g \in G_\mathbb{C}$. Then $H_{x\mathbb{C}} = gH_{w\mathbb{C}}g^{-1}$, and $g^{-1}y$ is fixed by $H_{w\mathbb{C}}$. So we may assume $x = w$. By (3.1), $V$ contains the trivial representation of $\mathfrak{h}$ precisely once. Therefore, the set of fixed points of $H_{w\mathbb{C}}$ is of dimension one. This proves (1).

Next consider (2). Let $x = gw$ with $g = (g_1, g_2) \in G_\mathbb{C}$ again. Let $S \subset \wedge^2 V_{1\mathbb{C}}$ be the three dimensional subspace corresponding to $y$. Then $g^{-1}S = g_1^{-1}S$ is fixed by $g_1^{-1}H_{x1\mathbb{C}}g_1 = H_{w1\mathbb{C}}$. Note that $G_2$ acts trivially on $\mathrm{Gr}(3, 10)$. So $g_1^{-1}S$ is $H_{w\mathbb{C}}$–invariant also. Therefore, $g_1^{-1}S$ is a sub $\mathfrak{h}_\mathbb{C}$–module. However, $\wedge^2 V_{1\mathbb{C}} \cong 6\Lambda \oplus 2\Lambda$. So there is only one three dimensional sub $\mathfrak{h}_\mathbb{C}$–module. Since $S_w$ also satisfies this condition, $y = g_1[S_w] = \Phi_3(x)$.

Statements (3), (4) are well known and were used in the proof of the Oppenheim conjecture for quadratic forms. $\square$

In the following theorem, let $x \in G_\mathbb{R} w$. Then $H_{x\mathbb{R}+}^0$ is generated by unipotent elements. Let $H_{x\mathbb{R}} \subset F \subset H_\mathbb{R}$ be the closed connected subgroup such that $\overline{H_{x\mathbb{R}} H_\mathbb{Z}} = FH_\mathbb{Z}$. By Ratner's theorem, such $F$ exists.

**Theorem (6.2)** (1) If $\Phi_2(x) \notin \mathbb{P}(\mathrm{Sym}^2 V_1^*)_\mathbb{Q}$ and $\Phi_3(x) \notin \mathrm{Gr}(3, 10)_\mathbb{Q}$, $F = \mathrm{SL}(5)_\mathbb{R} \times H_{x2\mathbb{R}}$ or $F = \mathrm{SL}(5)_\mathbb{R} \times \mathrm{SL}(3)_\mathbb{R}$.
(2) If $\Phi_1(x) \notin \mathbb{P}(\mathrm{Sym}^2 V_2^*)_\mathbb{Q}$, $\Phi_2(x) \notin \mathbb{P}(\mathrm{Sym}^2 V_1^*)_\mathbb{Q}$, and $\Phi_3(x) \notin \mathrm{Gr}(3, 10)_\mathbb{Q}$, $F = \mathrm{SL}(5)_\mathbb{R} \times \mathrm{SL}(3)_\mathbb{R}$.

*Proof.* Suppose $F = H_{x\mathbb{R}}$. Then $F$ is defined over $\mathbb{Q}$. Therefore, for any $\sigma \in \mathrm{Aut}(\mathbb{C}/\mathbb{Q})$, $H_{x\mathbb{C}}^\sigma = H_{x\mathbb{C}}$. Since $H_{x\mathbb{C}}^\sigma = H_{x^\sigma \mathbb{C}}$, $x^\sigma$ is fixed by $H_{x\mathbb{C}}$. So $x^\sigma$ is a scalar multiple of $x$ by Lemma (6.1). Since this is the case for all $\sigma$, $[x] \in \mathbb{P}(V)_\mathbb{Q}$. Since $\Phi_1, \Phi_2, \Phi_3$ are defined over $\mathbb{Q}$, $\Phi_1(x), \Phi_2(x), \Phi_3(x)$ are all $\mathbb{Q}$–rational points.

Suppose $F = H_{x1\mathbb{R}} \times H_{x2\mathbb{R}}$, or $H_{x1\mathbb{R}} \times \mathrm{SL}(3)_\mathbb{R}$. We show that $\Phi_3(x) \in \mathrm{Gr}(3, 10)_\mathbb{Q}$. Since the argument is similar, we only consider the first case.

For any $\sigma \in \mathrm{Aut}(\mathbb{C}/\mathbb{Q})$,

$$H^\sigma_{x1\mathbb{C}} \times H^\sigma_{x2\mathbb{C}} = H_{x^\sigma 1\mathbb{C}} \times H_{x^\sigma 2\mathbb{C}} = H_{x1\mathbb{C}} \times H_{x2\mathbb{C}}.$$

Since $G_2$ acts trivially on $\mathrm{Gr}(3, 10)$, $\Phi_3(x^\sigma)$ is fixed by $H_{x1\mathbb{C}}$. This implies $\Phi_3(x^\sigma) = \Phi_3(x)$. Since $\Phi_3$ is defined over $\mathbb{Q}$, $\Phi_3(x^\sigma) = \Phi_3(x)^\sigma = \Phi_3(x)$ by Lemma (6.1). Therefore, $\Phi_3(x) \in \mathrm{Gr}(3, 10)_\mathbb{Q}$.

By a similar argument, if $F = H_{x1\mathbb{R}} \times H_{x2\mathbb{R}}, H_{3\mathbb{R}} \times H_{x2\mathbb{R}}$, or $\mathrm{SL}(5)_\mathbb{R} \times H_{x2\mathbb{R}}$, $[\Phi_1(x)] \in \mathbb{P}(\mathrm{Sym}^2 V_2^*)_\mathbb{Q}$. Also if $F = H_{x3\mathbb{R}} \times H_{x2\mathbb{R}}$ or $H_{x3\mathbb{R}} \times \mathrm{SL}(3)_\mathbb{R}$, $[\Phi_2(x)] \in \mathbb{P}(\mathrm{Sym}^2 V_1^*)_\mathbb{Q}$.

By these considerations, conditions in (1), (2) force $F$ to become the given subgroups. $\square$

**Lemma (6.3)** *Let $x \in G_\mathbb{R} w$. Then for any non-zero real number $r$, there exists $h \in H_\mathbb{R}$ and a primitive integer point $a \in V_{1\mathbb{Z}}$ such that $F_{h^{-1}x}(a) = r$.*

*Proof.* We may assume $x = \lambda w$ where $\lambda \in \mathbb{R} \setminus \{0\}$. Since $F_{\lambda h^{-1} w}(a) = \lambda^{60} F_w(ha)$, the above condition is equivalent to $F_w(ha) = \lambda^{-60} r$. Put $t = -\lambda^{-20}(\frac{r}{2})^{\frac{1}{3}}$. Then

$$h = \left( \begin{pmatrix} t^{-1} & 0 & 0 & 0 & 0 \\ 0 & 1 & 0 & 0 & 0 \\ 0 & 0 & t & 0 & 0 \\ 0 & 0 & 0 & 1 & 0 \\ 0 & 0 & 0 & 0 & 1 \end{pmatrix}, I_3 \right)$$

and $a = {}^t\begin{pmatrix} 0 & 0 & 1 & 0 & 0 \end{pmatrix}$ satisfy the condition. $\square$

In the following theorem, $x \in G_\mathbb{R} w$.

**Theorem (6.4)** *If $[\Phi_2(x)] \notin \mathbb{P}(\mathrm{Sym}^2 V_1^*)_\mathbb{Q}$, and $\Phi_3(x) \notin \mathrm{Gr}(3, 10)_\mathbb{Q}$, the set of values of the cubic form $F_x(a)$ at primitive integer points is dense in $\mathbb{R}$.*

*Proof.* Let $r$ be a non-zero real number. We choose $h = (h', h'') \in H_\mathbb{R}$ and $a \in V_{1\mathbb{Z}}$ as in Lemma (6.3). By Theorem (6.2), there exist $h_1 = (h'_1, h''_1) \in H_{x\mathbb{R}}$ and $h_2 = (h'_2, h''_2) \in H_\mathbb{Z}$ such that $h'_1 h'_2$ is close to $h'$. Then

$$F_x(h'_2 a) = F_{h_2^{-1} x}(a) = F_{h_2^{-1} h_1^{-1} x}(a) = F_{(h'_2{}^{-1} h'_1{}^{-1}, 1)x}(a)$$

is close to

$$F_{(h'^{-1}, 1)x}(a) = F_{h^{-1}x}(a) = r.$$

Note that $F_{h^{-1}x}$ does not depend on the second component of $h$.

Since $h'_2 a \in V_{1\mathbb{Z}}$ is primitive, this proves the theorem. $\square$

Note that if $x = gw$ with $g = (g_1, g_2) \in G_\mathbb{R}$,

$$\Phi_3(x) = g_1 \Phi_3(w), \ [\Phi_2(x)] = g_1[\Phi_2(w)], \ F_x(a) = (\det g_1)^{15} (\det g_2)^3 F_w(g_1^{-1} a).$$

Therefore, writing down $F_w$, etc. explicitly, we get the statement of Theorem (0.2).

Akihiko Yukie  
Oklahoma State University  
Mathematics Department  
401 Math Science  
Stillwater OK 74078-1058 USA  
yukie@math.okstate.edu  
http://www.math.okstate.edu/~yukie